\documentclass{article}

\usepackage{amssymb} 
\usepackage{amsmath} 
\usepackage{latexsym} 
\usepackage{theorem} 
\usepackage{color}
\usepackage{epsfig}

\theorembodyfont{\itshape} 

\newtheorem{theorem}{Theorem}[section]

\newtheorem{definition}[theorem]{Definition}

\newtheorem{remark}[theorem]{Remark}
\newenvironment{proof}{{\par\addvspace{0.1cm}\noindent \bf Proof. }}{\hfill$\Box$\par\medskip} 
 
\def\n{m}

\def\an{{\alpha-\n}}
\def\a{\alpha}

\def\Om{\varOmega}
\def\RR{\mathbb{R}}

\title{\bf The incenter of a triangle as a cone isoperimetric center}
\author{Jun O'Hara}

\numberwithin{equation}{section}

\begin{document}

\maketitle

\begin{abstract} 
We show that the the image of the regular projection of a vertex of a cone over a triangle that minimizes the ratio of the cube of the area of the boundary of the cone and the square of the volume of the cone coincides with the incenter. 
\end{abstract}



\section{Introduction}
%
Recently Katsuyuki Shibata introduced a new kind of center of triagles, which he calls the {\em illuminating center} (\cite{Sh}). 
Speaking a concept, it is a point that maximizes the total brightness of a triangular park $\Om$ obtained by a light source on that point, namely, a point that 
maximizes $V_0(x)=\int_{\varOmega}|x-y|^{-2}\,d\mu(y)$, where $\mu$ is the standard Lesbegue measure of $\RR^2$. 
Unfortunately, $V_0(x)$ is not well-defined; it diverges for any point in $\Om$. 
In order to produce a well-defined potential, Shibata used the cut-off of the divergence of the integrand. 

In \cite{O1} the author introduced renormalization of the Riesz potential $\int_\Om{|x-y|}^{\an}\,d\mu(y)$ of a compact set $\Om$ in $\RR^\n$ (which is a closure of an open set) for $\a\le0$ to obtain a one-parameter family of {\em (renormalized) potentials} $V^{(\a)}_\Om$, and studied 
the points that attain the extremal values of $V^{(\a)}_\Om$, which we call the {\em $r^{\a-\n}$-centers} of $\Om$. 
The notion of ${r^{\an}}$-centers includes not only Shibata's illuminating center of a planar domain as an $r^0$-center, but also 
the center of mass of any comapct set $\Om\subset\RR^\n$ as $r^2$-center. This is because the center of mass $x_G$ is given by 
$x_G={\int_{\Om}y\,d\mu(y)}/{\int_{\Om}1\,d\mu(y)}\,,$ 
or equivalently by 
$\int_{\Om}(x_G-y)\,d\mu(y)=0\,,$
which implies that it can be characterized as a unique critical point of the map $V_\Om^{(\n+2)}:\RR^\n\ni x\mapsto\int_{\Om}{|x-y|}^2\,d\mu(y)\in\RR$. 

Shibata announced\footnote{at 2010 Autumn Meetings of the Mathematical Society of Japan} a theorem that an $r^{\a'}$-center of a non-obtuse triangle approaches the circumcenter as $\a'$ goes to $+\infty$ and to the incenter as $\a'$ goes to $-\infty$. 
The proof with more generality is given in \cite{O1}. 
%
Thus, we can give interpretations of the barycenters, the circumcenters, and the incenters of triangles as points that optimize a kind of potential or the limits of them. 

The motivation of the theorem in this note comes from the same philosophy; to express a center as a point that optimizes a kind of potential. We only deal with triangles in this note. For our potential, we use the ratio of the volume of the cone over $\Om$ and that of its boundary, with the former being squared and the latter cubed to make the ratio scale invariant. Then, the image of the regular projection of a vertex of a cone that optimizes this ratio is nothing but the incenter. 

\section{Cone isoperimetric center}
%
Let us start with general setting. 
Let $\Om$ be a compact set which is a closure of an open subset of $\RR^\n$ with a piecewise $C^1$ boundary $\partial\Om$. 
We assume that $\RR^\n$ is embedded in $\RR^{\n+1}$ in a standard way; $\RR^\n=\{(x_1,\ldots,x_\n,0)\in\RR^{\n+1}\,|\,x_i\in\RR\}$. 
Let $\Pi_h$ denote a level hyperplane in $\RR^{\n+1}$ with height $h>0$; $\Pi_h=\{x_{\n+1}=h\}$, and $C_p$ a cone over $\Om$ with vertex $p\in\Pi_h$; $C_p=\{tx+(1-t)p\,|\,x\in\Om, 0\le t\le 1\}$. 
Let $\pi:\RR^{\n+1}\to\RR^\n$ be the regular projection. 

\begin{definition} \rm 
\begin{enumerate}
\item Let $p_h$ be a point in $\Pi_h$ that attains the minimum value of a function $\Pi_h\ni p\mapsto \textrm{Vol}\,(\partial C_p)$. 
We call $\pi(p_h)$ a {\em cone isoperimetric center of $\Om$ of height $h$}. 
\item Let $p$ be a point in $\RR^{\n+1}_+=\{x_{\n+1}>0\}$ that attains the minimum value of a function 
$$
f(p)=\frac{\left(\textrm{Vol}\,(\partial C_p)\right)^{\n+1}}{\left(\textrm{Vol}\,(C_p)\right)^\n}.
$$
We call $C_p$ an {\em isoperimetrically optimal cone} and $\pi(p)$ a {\em cone isoperimetric center of $\Om$}. 
\end{enumerate}
\end{definition}

\begin{theorem}
Let $\Om$ be a triangle. 
\begin{enumerate}
\item The cone isoperimetric center of height $h$ coincides with the incenter for any $h>0$. 
\item The height of the isoperimetrically optimal cone is $2\sqrt2$ times the radius of the inscribed circle. 
\end{enumerate}
\end{theorem}
\begin{proof}
(1) Let $A,B$, and $C$ be vertices of the triangle, $S$ the area, $a,b$, and $c$ the lengths of the edges $BC, CA$, and $AB$ respectively. 
Let $P\in \Pi_h$ be a point and $D=\pi(P)$. 
Let $u,v$, and $w$ be the distances with signs between $D$ and the lines $\overline{BC}, \overline{CA}$, and $\overline{AB}$ respectively. 
The signs of $u,v$, and $w$ are given as follows. 
We put $u>0$ if $D$ and $A$ are in the same half-plane cut out by the line $\overline{BC}$. 
Then the area of the triangle $\Delta ABC$ is given by $S=\frac12(au+bv+cw)$, and the area of the boundary of the cone is given by 
$$
\textrm{Vol}\,(\partial C_P)=S+\frac12\left(a\sqrt{u^2+h^2}+b\sqrt{v^2+h^2}+c\sqrt{w^2+h^2}\,\right).
$$
Remark that the position of $D$ is determined uniquely by $u$ and $v$. 

It is obvious that a cone isoperimetric center of $\Delta ABC$ of height $h$ exists as $\textrm{Vol}\,(\partial C_P)$ goes to $+\infty$ as $|P|$ goes to $+\infty$. 
Let $D_h$ be a cone isoperimetric center of $\Delta ABC$ of height $h$, and $u_h,v_h$, and $w_h$ be the signed disntances between $D_h$ and the lines $\overline{BC}, \overline{CA}$, and $\overline{AB}$ respectively. 
Then the pair $(u_h,v_h)$ minimizes a function 
$$
F(u,v)=a\sqrt{u^2+h^2}+b\sqrt{v^2+h^2}+c\sqrt{\left(\frac{2S-au-bv}{c}\right)^2+h^2}. 
$$
\setlength\arraycolsep{1pt}
Therefore, when $(u,v,w)=(u_h,v_h,w_h)$ we have 
\[\begin{array}{rcl}
F_u(u,v)&=&\displaystyle \frac{au}{\sqrt{u^2+h^2}}+\frac{cw}{\sqrt{w^2+h^2}}\cdot\left(-\frac{a}{c}\right)=0,\\[4mm]
F_v(u,v)&=&\displaystyle \frac{bv}{\sqrt{v^2+h^2}}+\frac{cw}{\sqrt{w^2+h^2}}\cdot\left(-\frac{b}{c}\right)=0,
\end{array}\]
which implies 
\begin{equation}\label{f_equi_angle}
\frac{u}{\sqrt{u^2+h^2}}=\frac{v}{\sqrt{v^2+h^2}}=\frac{w}{\sqrt{w^2+h^2}}. 
\end{equation}

Remark that the above holds only when $u,v$, and $w$ are all positive, implying that $D_h$ is in the interior of $\Delta ABC$. 
The equation \eqref{f_equi_angle} means that three angles between the $xy$ plane and three planes through $PAB$, $PBC$, and $PCA$ are all equal. 
Therefore, for each pair of the three planes above mentioned, there is the symmetry in a plane orthogonal to the $xy$ plane that contains the intersection of the pair. 
Looking from above, you can see that three lines $D_hA$, $D_hB$, and $D_hC$ 
are the angle bisectors of $\angle A$, $\angle B$, and $\angle C$ respectively, which means that $D_h$ is the incenter of $\Delta ABC$. 

\medskip
(2) The statement follows from elementary calculus. 
Let $r$ be the radius of the inscribed circle. 
Put $P_h=\pi^{-1}(D_h)\cap \Pi_h$, then 
$$
\textrm{Vol}\,(\partial C_{P_h})=\frac12\left((a+b+c)r+(a+b+c)\sqrt{r^2+h^2}\right)
=S\left(1+\sqrt{1+\left(\frac hr\right)^2}\,\right).
$$
As $\textrm{Vol}\,(C_{P_h})=\frac13Sh$, 
$$
f(P_h)=\frac{\left(\textrm{Vol}\,(\partial C_{P_h})\right)^3}{\left(\textrm{Vol}\,(C_{P_h})\right)^2}
=9S\frac{\left(1+\sqrt{1+\left(\frac hr\right)^2}\,\right)^3}{h^2}
=\frac{9S}{r^2}\cdot\frac{\left(1+\sqrt{1+\left(\frac hr\right)^2}\,\right)^3}{\left(\frac hr\right)^2}. 
$$
Since $\varphi(t)=\frac{\left(1+\sqrt{1+t^2}\,\right)^3}{t^2}$ $(t>0)$ takes the minimum at $t=2\sqrt2$, it completes the proof. 
\end{proof}

\begin{remark}\rm 
As an example below shows, the cone isoperimetric center of height $h$ is not identically the same, and the cone isoperimetric center does not coincide with the limit of $r^{\a'}$-center as $\a'$ goes to $-\infty$ in general. 

Let us call a point an {\em asymptotic $r^{-\infty}$-center} of $\Om$ if it is the limit of a convergent sequence of $r^{\a_i'}$-centers with $\a_i'\to-\infty$ as $i\to+\infty$. 
We showed in \cite{O1} that an asymptotic $r^{-\infty}$-center is a {\em max-min point} of $\Om$, by which we mean a point that attains the supremum of a map $\mathbb R^\n\ni x\mapsto \min_{y\in\overline{{\Om}^c}}|y-x|\in\RR$, where $\overline{{\Om}^c}$ denotes the closure of the complement of $\Om$. 
We remark that an $r^\an$-center $(\a<0)$ and a max-min point are not necessarily unique. 
To see this, it is enough to consider a disjoint union of two rectangles, say, $\Om'=\{(\xi,\eta)\,|\,1\le|\xi|\le2, \,|\eta|\le2\}$. 

Let $\Om$ be a trapezoid given by $\Om=\{(\xi,\eta)\,|\,0\le\xi\le2, \,|\eta|\le1+\frac12\xi\}$. 
It is easy to see that a cone isoperimetric center of height $h$ is on the $\xi$-axis for any $h$. 
Let $(\xi_h,0)$ be the coordinates of it. 
Numerical experiment shows that $\xi_1\sim0.9169, \xi_2\sim0.9079, \xi_3\sim0.9045$, and $\xi_4\sim0.9031$, and the minimum of the ratio $f$ is attained at $h\sim3.250$ when $\xi_h\sim0.90405$. 
On the other hand, an asymptotic $r^{-\infty}$-center is $(1,0)$. This is because the set of max-min points is $\{(1,\eta)\,|\,|\eta|\le\frac32-\frac{\sqrt5}2\}$ whereas any $r^{\a'}$-center is contained in $\{(\xi,0)\,|\,1\le\xi\le\frac74\}$ for any $\a'$ by the symmetry argument (based on the moving plane method \cite{GNN}) explained in \cite{O1}, and the point $(1,0)$ is the unique intersection point of these sets. 
\end{remark}
%

Department of Mathematics and Information Sciences, 

Tokyo Metropolitan University, 

1-1 Minami-Ohsawa, Hachiouji-Shi, Tokyo 192-0397, JAPAN. 

E-mail: ohara@tmu.ac.jp


\begin{thebibliography}{O} 
\bibitem[GNN]{GNN}B.~Gidas, W. M. ~Ni, and L.~Nirenberg, {\em Symmetry and related properties via the maximum principle,} Comm. Math. Phys. {\bf 68} (1979), 209\,--\,243.

\bibitem[O]{O1}J.~O'Hara, {\em Renormalization of $r^{\bullet}$-potentials and generalization of dual volumes and centers,} arXiv:1008.2731. 


\bibitem[S]{Sh}K.~Shibata, {\em Where should a streetlight be placed in a triangle-shaped park?
Elementary integro-differential geometric optics,} available at 
http://www1.rsp.fukuoka-u.ac.jp/kototoi/shibataaleph-sjs.pdf
\end{thebibliography}
\end{document}